\documentclass[10pt]{article}
\usepackage[a4paper,left=2cm,right=2cm,top=2.5cm,bottom=2.5cm]{geometry}
\usepackage{amsmath,amssymb}
\usepackage{hyperref}
\usepackage{setspace}
\usepackage[utf8]{inputenc}
\usepackage[english]{babel}
\usepackage{booktabs,multirow}
\usepackage{blkarray}
\usepackage{lscape}
\usepackage[T1]{fontenc}
\usepackage{authblk,color}
\usepackage[english]{babel}
\usepackage{amsfonts, amstext, amsthm, booktabs, dcolumn}
\usepackage{graphicx}
\usepackage{float}
\usepackage{caption}
\usepackage{centernot}
\usepackage{subcaption}

\definecolor{gr}{rgb}{0.7, 0.0, 0.49}

\newtheorem{theorem}{\bf Theorem}[section]

\newcommand\floor[1]{\left\lfloor#1\right\rfloor}

\begin{document}
\title{Generalizations of Distributions Related to ($k_1,k_2$)-runs}
%\title{Distributions Related to Generalized $(k_1,k_2)$-runs}
%\title{On Distribution of Generalized $(k_1,k_2)$-runs}
\author[*]{A. N. Kumar}
\author[**]{N. S. Upadhye}
\affil[ ]{\small Department of Mathematics, Indian Institute of Technology Madras,}
\affil[ ]{\small Chennai-600036, India.}
\affil[*]{\small Email: amit.kumar2703@gmail.com}
\affil[**]{\small Email: neelesh@iitm.ac.in}
\date{}
\maketitle

\begin{abstract}
\noindent
The paper deals with three generalized dependent setups arising from a sequence of Bernoulli trials. Various distributional properties, such as probability generating function, probability mass function and moments are discussed for these setups and their waiting time. Also, explicit forms of probability generating function and probability mass function are obtained. Finally, two applications to demonstrate the relevance of the results are given.
\end{abstract}

\noindent

\section{Introduction}
Runs and patterns play a crucial role in applied statistics and have numerous applications, for example, reliability theory (see Fu \cite{FU} and Fu and Hu \cite{FH}), nonparametric hypothesis testing (Balakrishnan and Koutras \cite{BK}), DNA sequence analysis (Fu {\em et al.} \cite{FLC}), statistical testing (Balakrishnan {\em et al.} \cite{BMA}), computer science (Sinha {\em et al.} \cite{SSD}), the start-up demonstration theory (Balakrishnan {\em et al.} \cite{BKM1, BKM2}) and quality control (Moore \cite{moore}) among many others.\\
A run can be defined as an occurrence of specific patterns of failures or successes or both in a sequence of Bernoulli trials. In particular, a pattern of consecutive successes of length $k$ is considered by Philippou {\em et al.} \cite{PGP} and described geometric and negative binomial distribution of order $k$. Also, Philippou and Makri \cite{PM} discussed binomial distribution of order $k$. Later, Huang and Tsai \cite{HT} extended the pattern by observing at least $k_1$  consecutive failures followed by at least $k_2$ consecutive successes and studied a modified binomial distribution of order $k$ or $(k_1,k_2)$-runs. Recently, Dafnis {\em et al.} \cite{DAP} also considered three types of $(k_1,k_2)$-runs which include the pattern discussed in Huang and Tsai \cite{HT}. Though there have been several studies on this topic, still there are many problems which can not be studied based on the available literature. For example, (i) let us consider the quality control problem in which the system is said to be in control, whenever, (on the control chart) not more than two consecutive points exceed the control limits and at least three succeeding points are inside the control limits (see {\bf (T1)} below with $\ell_1=1$, $k_1=2$ and $\ell_2=3$). Similarly, (ii) consider a climatology problem, in which, climatologist is interested in knowing the distribution of at least two consecutive rainy days followed by exactly five consecutive dry days (see {\bf (T2)} below with $\ell_1=2$ and $\ell_2=k_2=5$). Also, there are several such problems that occur in brand switching, learning, reliability and queuing models. Hence, there is a need to generalize the results related to ($k_1,k_2$)-runs.\\
In this paper, we generalize $(k_1,k_2)$-runs to include the following patterns, for $1 \le \ell_1\le k_1$ and $1 \le \ell_2\le k_2$,
\begin{itemize}
\item[{\bf (T1)}] at least $\ell_1$ and at most $k_1$ consecutive 0's followed by at least $\ell_2$ consecutive 1's.
\item[{\bf (T2)}] at least $\ell_1$ consecutive 0's followed by at least $\ell_2$ and at most $k_2$ consecutive 1's.
\item[{\bf (T3)}] at least $\ell_1$ and at most $k_1$ consecutive 0's followed by at least $\ell_2$ and at most $k_2$ consecutive 1's.
\end{itemize}
Note that {(\bf T1)}, {(\bf T2)} and {(\bf T3)} contain various ($k_1,k_2$)-runs. For example,
\begin{enumerate}
\item if $\ell_1=k_1$ then {(\bf T1)} leads to, exactly $\ell_1$ consecutive 0's followed by at least $\ell_2$ consecutive 1's,
\item if $\ell_2=k_2$ then {(\bf T2)} leads to, at least $\ell_1$ consecutive 0's followed by exactly $\ell_2$ consecutive 1's,
\item if $\ell_1=\ell_2=1$ then {(\bf T3)} leads to, at most $k_1$ consecutive 0's followed by at most $k_2$ consecutive 1's,
\item if $\ell_1=k_1$ and $\ell_2=k_2$ then {(\bf T3)} leads to, exactly $k_1$ consecutive 0's followed by exactly $k_2$ consecutive 1's
\end{enumerate}
and similarly, other special cases can be seen by choosing the values for $\ell_1$ and $\ell_2$, $k_1$ and $k_2$ appropriately. Dafnis {\em et al.} \cite{DAP} considered two special cases of ${\bf (T3)}$, namely, (i) $\ell_1=1=\ell_2$ and (ii) $\ell_1=k_1$ and $\ell_2=k_2$. \\
Next, let $\zeta_1,\zeta_2,\dotsc,\zeta_n$ be a finite sequence of independent Bernoulli trials with success (denoted by 1) probability $p$ and failure (denoted by 0) probability $q=1-p$. Then, define
\begin{align*}
I_{s}^{(m)}\hspace{-0.1cm}&:= \hspace{-0.1cm}\left\{\hspace{-0.15cm}
\begin{array}{ll}
(1\hspace{-0.08cm}-\hspace{-0.08cm}\zeta_{m})\dotsb(1\hspace{-0.08cm}-\hspace{-0.08cm}\zeta_{m+\ell_1-1})(1\hspace{-0.08cm}-\hspace{-0.08cm}\zeta_{m+\ell_1})\dotsb(1\hspace{-0.08cm}-\hspace{-0.08cm}\zeta_{m+s+\ell_1-1})\zeta_{m+s+\ell_1}\dotsb \zeta_{m+s+\ell_1+\ell_2-1}, & \hspace{-0.1cm}m = 1,\\
\zeta_m(1\hspace{-0.08cm}-\hspace{-0.08cm}\zeta_{m+1})\dotsb(1\hspace{-0.08cm}-\hspace{-0.08cm}\zeta_{m+\ell_1})(1\hspace{-0.08cm}-\hspace{-0.08cm}\zeta_{m+\ell_1+1})\dotsb(1\hspace{-0.08cm}-\hspace{-0.08cm}\zeta_{m+s+\ell_1})\zeta_{m+s+\ell_1+1}\dotsb \zeta_{m+s+\ell_1+\ell_2},&\hspace{-0.1cm} 2\hspace{-0.07cm}\le\hspace{-0.07cm} m \hspace{-0.07cm}\le \hspace{-0.07cm}n\hspace{-0.07cm}-\hspace{-0.07cm}\ell_1\hspace{-0.07cm}-\hspace{-0.07cm}\ell_2,
\end{array}\right.\\
J_{t}^{(m)}\hspace{-0.12cm} &:=\hspace{-0.12cm} (\hspace{-0.02cm}1\hspace{-0.08cm}-\hspace{-0.08cm}\zeta_{m}\hspace{-0.02cm})\dotsb(\hspace{-0.02cm}1\hspace{-0.08cm}-\hspace{-0.08cm}\zeta_{m+\ell_1-1}\hspace{-0.02cm})\zeta_{m+\ell_1}\dotsb\zeta_{m+\ell_1+\ell_2-1}\zeta_{m+\ell_1+\ell_2}\dotsb\zeta_{m+t+\ell_1+\ell_2-1}(\hspace{-0.02cm}1\hspace{-0.08cm}-\hspace{-0.08cm}\zeta_{m+t+\ell_1+\ell_2}\hspace{-0.02cm}), 1\hspace{-0.07cm}\le\hspace{-0.07cm} m \hspace{-0.07cm}\le\hspace{-0.07cm} n\hspace{-0.07cm}-\hspace{-0.1cm}\ell_1\hspace{-0.07cm}-\hspace{-0.07cm}\ell_2,\\
K_{s,t}^{(m)}&:=\left\{
\begin{array}{ll}
(1-\zeta_{m})\dotsb(1-\zeta_{m+\ell_1-1})(1-\zeta_{m+\ell_1})\dotsb(1-\zeta_{m+s+\ell_1-1})\zeta_{m+s+\ell_1} & \\
~~~~~~~~~~\dotsb \zeta_{m+s+\ell_1+\ell_2-1}\zeta_{m+s+\ell_1+\ell_2}\dotsb\zeta_{m+s+t+\ell_1+\ell_2-1}(1-\zeta_{m+s+t+\ell_1+\ell_2}),& m = 1,\\
\zeta_m(1-\zeta_{m+1})\dotsb(1-\zeta_{m+\ell_1})(1-\zeta_{m+\ell_1+1})\dotsb(1-\zeta_{m+s+\ell_1})\zeta_{m+s+\ell_1+1} & \\
~~~~~~~~~~\dotsb \zeta_{m+s+\ell_1+\ell_2}\zeta_{m+s+\ell_1+\ell_2+1}\dotsb\zeta_{m+s+t+\ell_1+\ell_2}(1-\zeta_{m+s+t+\ell_1+\ell_2+1}),& 2\hspace{-0.05cm}\le \hspace{-0.05cm}m\hspace{-0.05cm} \le\hspace{-0.05cm} n\hspace{-0.05cm}-\hspace{-0.05cm}\ell_1\hspace{-0.05cm}-\hspace{-0.05cm}\ell_2\hspace{-0.05cm}-\hspace{-0.05cm}1,
\end{array} \right.\\
I_m:&=\max_{0\le s\le k_1-\ell_1}I_{s}^{(m)},\quad\quad\quad\quad J_m:=\max_{0\le t \le k_2-\ell_2}J_{t}^{(m)},\quad\quad\quad\quad K_m:=\max_{\substack{0\le s\le k_1-\ell_1\\0\le t \le k_2-\ell_2}}K_{s,t}^{(m)}.
\end{align*}
Note that $I_s^{(m)}$, $J_t^{(m)}$ and $K_{s,t}^{(m)}$ denote a $(k_1,k_2)$-run of type {\bf (T1)}, {\bf (T2)} and {\bf (T3)}, respectively. For example, if $\ell_1=1$, $k_1=3$, $\ell_2=2$, $k_2=2$ and $m=1$ then the possible patterns for {\bf (T1)}, {\bf (T2)} and {\bf (T3)} are $\{\bf 011,~0011,~00011\}$, $\{\bf 0110\}$ and $\{\bf 0110,~00110,~000110\}$, respectively. This can also be verified using the definitions of $I_s^{(1)}$, $J_t^{(1)}$ and $K_{s,t}^{(1)}$.
 %Observe that we do not require a success before the pattern for {\bf (T1)} and {\bf (T3)} as $m=1$ otherwise it is needed. Now, from the definition of $I_s^{(m)}$, $J_t^{(m)}$ and $K_{s,t}^{(m)}$, we have
%\begin{align*}
%I_s^{(1)}&=(1-\zeta_1)\dotsb(1-\zeta_{s+1})\zeta_{s+2}\zeta_{s+3},\quad J_{t}^{(1)}=(1-\zeta_1)\zeta_2\dotsb\zeta_{t+3}(1-\zeta_{t+4})\\
%{\rm and}~~~~~~~~~~~~~~~~~~K_{s,t}^{(1)}&=(1-\zeta_1)\dotsb(1-\zeta_{s+1})\zeta_{s+2}\dotsb\zeta_{s+t+3}(1-\zeta_{s+t+4}).~~~~~~~~~~~~~~~~~~~~~~~~~~~~~~~~~~~~~~~~~~~~~~~~~~~~~~~~~~~~~~~
%\end{align*}
%Now, note that if $I_s^{(1)}=J_t^{(1)}=K_{s,t}^{(1)}=1$ for the range $s=0$ to $k_1-\ell_1=2$ and $t=0=k_2-\ell_2$ then it can seen that $\{\bf 011,~0011,~00011\}$, $\{\bf 0110\}$ and $\{\bf 0110,~00110,~000110\}$ occur for $I_s^{(m)}$, $J_t^{(m)}$ and $K_{s,t}^{(m)}$, respectively, as expected. 
\\
Next, let $H_{\ell_1,k_1,\ell_2}^n$, $H_{\ell_1,\ell_2,k_2}^n$ and $H_{\ell_1,k_1,\ell_2,k_2}^n$ be the number of occurrences for {\bf (T1)}, {\bf (T2)} and {\bf (T3)} type events, respectively. Then, random variable representation of $H_{\ell_1,k_1,\ell_2}^n$, $H_{\ell_1,\ell_2,k_2}^n$ and $H_{\ell_1,k_1,\ell_2,k_2}^n$ can be seen as follows:
\begin{align*}
H_{\ell_1,k_1,\ell_2}^n&=\sum_{m=1}^{n-\ell_1-\ell_2}I_m,\quad\quad H_{\ell_1,\ell_2,k_2}^n=\sum_{m=1}^{n-\ell_1-\ell_2}J_m~\quad\quad{\rm and}~\quad\quad H_{\ell_1,k_1,\ell_2,k_2}^n=\sum_{m=1}^{n-\ell_1-\ell_2-1}K_m.
\end{align*}
Now, let us consider a particular realization in a sequence of 20 Bernoulli trials given by
$${\bf 0~0~1~1~1~1~0~1~1~0~0~0~1~0~1~0~0~0~1~1}.$$
Here, note that
\begin{itemize}
\item[{\bf (T1)}]  $H_{1,1,1}^{20}=2$, $H_{1,2,2}^{20}=2$, $H_{2,2,3}^{20}=1$ and $H_{1,2,1}^{20}=3$.
\item[{\bf (T2)}] $H_{1,1,2}^{20}=3$, $H_{3,1,2}^{20}=1$, $H_{2,2,2}^{20}=0$ and $H_{1,4,4}^{20}=1$.
\item[{\bf (T3)}] $H_{1,1,1,1}^{20}=1$, $H_{1,2,2,2}^{20}=1$, $H_{1,1,1,2}^{20}=2$ and $H_{1,2,1,2}^{20}=2$.
\end{itemize}
%Next, runs and patterns play a crucial role in applied statistics and has numerous applications, for example, reliability theory (see Fu \cite{FU} and Fu and Hu \cite{FH}), nonparametric hypothesis testing (Balakrishnan and Koutras \cite{BK}), DNA sequence analysis (Fu {\em et al.} \cite{FLC}), statistical testing (Balakrishnan {\em et al.} \cite{BMA}), computer science (Sinha {\em et al.} \cite{SSD}) and quality control (Moore \cite{moore}) among many others.
For more details about runs and patterns, we refer the reader to Aki \cite{AKI}, Aki {\em et al.} \cite{AKH}, Antzoulakos {\em et al.} \cite{ABK}, Antzoulakos and Chadjiconstantinidis \cite{AC}, Balakrishnan and Koutras \cite{BK}, Dafnis {\em et al.} \cite{DAP}, Fu and Koutras \cite{FK}, Koutras \cite{KMV,koutras} and Makri {\em et al.} \cite{MPP} and references therein.\\
This paper is organized as follows. In Section \ref{four:DR}, we obtain the double probability generating function (PGF) and waiting time for $H_{\ell_1,k_1,\ell_2}^n$, $H_{\ell_1,\ell_2,k_2}^n$ and $H_{\ell_1,k_1,\ell_2,k_2}^n$. Using double PGF, we derive recursive relation in PGF, probability mass function (PMF) and moments and also derive an explicit form of PGF and PMF. Next, using double PGF for waiting time, we obtain the PGF, recursive relations in PMF and moments. Finally, we obtain the double PGF under Markov dependent trials. In Section \ref{APP}, we demonstrate the relevance of the results through some interesting applications. In Section \ref{four:CR}, we point out some relevant remarks.

\section{Distributions Related to  $H_{\ell_1,k_1,\ell_2}^n$, $H_{\ell_1,\ell_2,k_2}^n$ and $H_{\ell_1,k_1,\ell_2,k_2}^n$}\label{four:DR}
In this section, we discuss various distributional properties such as PGF, PMF and moments for $H_{\ell_1,k_1,\ell_2}^n$, $H_{\ell_1,\ell_2,k_2}^n$ and $H_{\ell_1,k_1,\ell_2,k_2}^n$ and their waiting time.\\
The method used can be formulated in the following way. Let $Y_n$ be a random variable related to ($k_1,k_2$)-runs. Then, we can define a Markov chain $\{Z_t,~t \ge 0\}$ on discrete space $\Omega$ (which can be partitioned into discrete subspaces $\{0,1,2,\dotsc,r\}$ of maximum length $\varepsilon_n$ and contains one and only one ($k_1,k_2$)-event) such that ($k_1,k_2$)-runs has occurred $v$ times if and only if Markov chain is in $v$-th discrete subspace (say $E_v=\{E_{v,0},E_{v,1},\dotsc,E_{v,r}\}$ such that $\Omega =\cup_{v\ge 0}E_v$). Now, assume $A$ and $B$ be $(r+1)\times(r+1)$ matrices when ($k_1,k_2$)-runs are observed from $v$ to $v$ and $v$  to $v+1$ times, respectively. Let $\phi_n(\cdot)$ and $\Phi(\cdot,\cdot)$ be the single and double generating function of $Y_n$ and $H_j(\cdot)$ and $H(\cdot,\cdot)$ be the single and double generating function of $j$-th waiting time for $Y_n$. Then, the double generating function for $Y_n$ and its waiting time is given by
\begin{equation}
\Phi(t,z)=\sum_{j=0}^{\infty}\phi_j(t)z^j=\kappa_0 (\vartheta(z,t))^{-1} {\bf 1}^t\label{four:dpgf}
\end{equation}
and
\begin{equation}
H(t,z)=\sum_{j=0}^{\infty}H_j(t)z^j=1+tz \kappa_0 (\vartheta(t,z))^{-1}B{\bf 1}^t\label{four:wtdpgf}
\end{equation}
respectively, where $\kappa_0$ is the initial distribution, $\vartheta(z,t)={\bf I}-z(A+tB)$ be $(r+1)\times(r+1)$ matrix, ${\bf 1}^t$ is the transpose of row matrix $(1,1,\dotsc,1)$ with $(r+1)$ entries and ${\bf I}$ is $(r+1)\times(r+1)$ identity matrix. For more details, we refer the reader to Antzoulakos {\em et al.} \cite{ABK} and Dafnis {\em et al.} \cite{DAP}.\\
Let us define some notations as
$$a(p):=q^{\ell_1}p^{\ell_2},~~\ell:=\ell_1+\ell_2,~~m_1:=k_1-\ell_1+1,~~m_2:=k_2-\ell_2+1,$$
$\rho_r$ is the $r$-th waiting time for $(k_1,k_2)$-runs, $p_{\cdot,n}$ and $g_r(\cdot)$ be the PMF of $(k_1,k_2)$-runs and $\rho_r$, respectively. Also, define $\mu_{n,j}$ and ${\tilde{\mu}}_{r,j}$ be the $j$-th (non-central) moment of $(k_1,k_2)$-runs and $\rho_r$, respectively, where $n$ denotes the number of Bernoulli trials.

\subsection{Distribution of $H_{\ell_1,k_1,\ell_2}^n$ and its Waiting Time}\label{four:D1}
Recall that $H_{\ell_1,k_1,\ell_2}^n$ is the number of occurrences of (at least $\ell_1$) at most $k_1$ consecutive 0's followed by at least $\ell_2$ consecutive 1's. Here, $r= k_1+\ell_2+1$ and $k_1^+$ is the element after $k_1$ consecutive 0's (if failures occur) in $\{0,1,\dotsc,k_1,k_1^+=k_1+1,k_1+2,\dotsc,k_1+\ell_1+1\}$. It is easy to see that ${\mathbb P}\big(H_{\ell_1,k_1,\ell_2}^0=0\big)=1$ and $\varepsilon_n:=\sup\big\{x:{\mathbb P}\big(H_{\ell_1,k_1,\ell_2}^n=x\big)>0\big\}=\floor{n/\ell}$. Therefore, $\kappa_0=(1,0,\dotsc,0)_{1\times (k_1+\ell_2+2)}$, $A=[a_{i,j}]_{(k_1+\ell_2+2)\times(k_1+\ell_2+2)}$ with non-zero entries
\begin{itemize}
\item $a_{i,1}=p$ and $a_{i,i+1}=q$ for $1 \le i \le \ell_1$,
\item $a_{i,k_1+3}=p$ and $a_{i,i+1}=q$ for $\ell_1+1 \le i \le k_1+1$,
\item $a_{k_1+2,1}=p$ and $a_{k_1+2,k_1+2}=q$,
\item $a_{i,2}=q$ for $k_1+3 \le i \le k_1+\ell_2+2$ and $a_{i,i+1}=p$ for $k_1+3 \le i \le k_1+\ell_2$,
\item $a_{k_1+\ell_2+2,k_1+\ell_2+2}=p$
\end{itemize}
and $B=[b_{i,j}]_{(k_1+\ell_2+2)\times(k_1+\ell_2+2)}$ is the matrix of non-zero entry $b_{k_1+\ell_2+1,k_1+\ell_2+2}=p$. Hence, using \eqref{four:dpgf}, it can be easily verified that
\vspace{-0.01cm}
\begin{equation}
\Phi(t,z)=\sum_{n=0}^{\infty}\phi_n(t)z^n=\frac{1}{1-z-(qz)^{\ell_1}(pz)^{\ell_2}(1-(qz)^{k_1-\ell_1+1})}= \frac{1}{1-z-a(p)z^\ell(t-1)\left(1-(qz)^{m_1}\right)}.\label{four:ddpgf2}
\end{equation}
Let us illustrate the result for $\ell_1=1$, $k_1=2$ and $\ell_2=2$. Here, we have $\kappa_0=(1,0,0,0,0,0)$, 
$$A=
\begin{pmatrix}
pe_1&qe_1&qe_2&qe_3&p(e_2+e_3)&{\bf 0}\\
pe_1&q(e_2+e_3)&{\bf 0}&qe_1&{\bf 0}&pe_3\\
\end{pmatrix}
~~~{\rm and}~~~B=
\begin{pmatrix}
{\bf 0}&pe_5\\
{\bf 0}&{\bf 0}\\
\end{pmatrix},
$$
where ${\bf 0}$ is a zero vector or matrix of appropriate length, $e_5=(0,0,0,0,1)^t$ and $e_k=(0,\dotsc,0,1,0,\dotsc,0)^t$, $k=1,2,3$ are $3\times 1$ column vectors. Hence, with some algebraic calculations, it can be verified that
$$\Phi(t,z)=\kappa_0 ({\bf I}-z(A+tB))^{-1} {\bf 1}^t=\frac{1}{1-z+(qz)(pz)^2(t-1)(1-(qz)^2)}.$$
This expression is same as \eqref{four:ddpgf2} for $\ell_1=1$, $k_1=2$ and $\ell_2=2$, as expected.\\ 
Next, using \eqref{four:ddpgf2}, we have the following results.
\begin{theorem}\label{four:thh2}
The recursive relation in PGF, PMF and moments of $H_{\ell_1,k_1,\ell_2}^n$, for $n \ge \ell$, are given by
\begin{itemize}
\item[(i)] $\phi_n(t)=\phi_{n-1}(t)+a(p)(t-1)\left[\phi_{n-\ell}(t)-q^{m_1}\phi_{n-\ell-m_1}(t)\right]$\\
with initial condition $\phi_n(t)=1,$ for $n \le \ell-1$.
\item[(ii)] $p_{m,n}=p_{m,n-1}+a(p)\left[p_{m-1,n-\ell}-p_{m,n-\ell}-q^{m_1}(p_{m-1,n-\ell-m_1}-p_{m,n-\ell-m_1})\right]$\\
with initial conditions $p_{0,n} =1$ and $p_{m,n}=0,~m>0$ for $n \le \ell-1$.
\item[(iii)] $\mu_{n,j}=\displaystyle{\mu_{n-1,j}+a(p)\sum_{k=0}^{j-1}{j \choose k}\left[\mu_{n-\ell,k}-q^{m_1}\mu_{n-\ell-m_1,k}\right]},~~~~$ for $j \ge 1$ \vspace{0.25cm}\\
with initial conditions $\mu_{n,0}=1$ and $\mu_{n,j}=0$ for all $j \ge 1$ and $n \le \ell-1$.
\end{itemize}
\end{theorem}
\noindent
{\bf Proof}. From \eqref{four:ddpgf2}, $(i)$ follows and using the definition of PGF, $(ii)$ follows. Substituting $t=e^x=\sum_{m=0}^{\infty}x^m/m!$ in $(i)$ and comparing the coefficient of $x^m/m!$, $(iii)$ follows. \qed\\
Next, we obtain an explicit form of PGF and PMF using Theorem \ref{four:thh2}.
\begin{theorem}\label{four:ex1}
Assume the conditions of Theorem \ref{four:thh2} hold, then PGF and PMF of $H_{\ell_1,k_1,\ell_2}^n$ are given by
\begin{itemize}
\item[(i)] $\phi_n(t) = \displaystyle{\sum_{u=0}^{\floor{\frac{n}{\ell}}}\sum_{v=0}^{\floor{\frac{n-u \ell}{\ell+m_1}}}{n-u (\ell-1)-v(\ell+m_1-1) \choose n-u \ell-v(\ell+m_1),u,v} (-1)^{v}q^{v m_1}(a(p)(t-1))^{u+v}}$.
\item[(ii)] $p_{m,n}=\displaystyle{\sum_{u=0}^{\floor{\frac{n}{\ell}}}\sum_{v=0}^{\floor{\frac{n-u \ell}{\ell+m_1}}}{n-u (\ell-1)-v(\ell+m_1-1) \choose n-u \ell-v(\ell+m_1),u,v}{u+v \choose m}(-1)^{u-m} q^{v m_1}a(p)^{u+v} },$
\end{itemize}
where ${n \choose u_1,u_2, \dotsc,u_s}=\frac{n!}{u_1!u_2!\dotsb u_s!}$.
\end{theorem}
\noindent
{\bf Proof}. $(i)$ For $(t,z) \in \big\{|t|\le1,~|z|<1 ~{\rm and}~|z+a(p)z^\ell(t-1)(1-(qz)^{m_1})|<1\big\}$, \eqref{four:ddpgf2} can be written as
\vspace{-0.2cm}
$$\Phi(t,z)=\sum_{n=0}^{\infty}\left(z+a(p)z^\ell(t-1)(1-(qz)^{m_1})\right)^n.$$
\vskip -2.5ex
\noindent
Now, using binomial expansion and interchanging summations, we get the required result.
\newpage
\noindent
$(ii)$ Following the steps similar to $(i)$ with recursive relation $(ii)$ of Theorem \ref{four:thh2}, the proof follows.\qed\\
Next, using \eqref{four:wtdpgf} with some algebraic manipulations, it can be easily verified that
\begin{equation}
H(t,z)=1+\sum_{r=1}^{\infty}\left(\frac{a(p)t^\ell(1-(qt)^{m_1})}{1-t+a(p)t^\ell(1-(qt)^{m_1})}\right)^r z^r.\label{four:wt2}
\end{equation}
Hence, using \eqref{four:wt2}, we have the following theorem.

\begin{theorem}\label{four:Wth1}
Let $\delta_{i,j}$ denote Kronecker delta function. The PGF, PMF and moments of $\rho_{r}$, for $r \ge 1$, are given by
\begin{itemize}
\item[(i)] $\begin{aligned}[t]
H_r(t)=\left(\frac{a(p)t^\ell(1-(qt)^{m_1})}{1-t+a(p)t^\ell(1-(qt)^{m_1})}\right)^r.
\end{aligned}$
\item[(ii)] $g_r(m)=g_r(m-1)+a(p)\left[g_{r-1}(m-\ell)-g_r(m-\ell)-q^{m_1}(g_{r-1}(m-\ell-m_1)-g_r(m-\ell-m_1))\right],$\\
for $m \ge \ell r$ with initial condition $g_0(m)=\delta_{m,0}$, $g_r(m)=0$ for $m \le \ell r-1$.
\item[(iii)] $\tilde{\mu}_{r,j}=\displaystyle{\sum_{k=0}^{j}{j \choose k}[\tilde{\mu}_{r,k}+a(p)(\ell^{j-k}-q^{m_1}(\ell+m_1)^{j-k})(\tilde{\mu}_{r-1,k}-\tilde{\mu}_{r,k})],\quad j \ge 1}$ \vspace{0.25cm}\\
with initial condition $\tilde{\mu}_{0,i}=\delta_{i,0}$.
\end{itemize}
\end{theorem}
\noindent
{\bf Proof}. Following the steps similar to the proof of Theorem \ref{four:thh2}, the results follow.\qed

\subsection{Distribution of $H_{\ell_1,\ell_2,k_2}^n$ and its Waiting Time}\label{four:D2}
Recall that $H_{\ell_1,\ell_2,k_2}^n$ is the number of occurrences of at least $\ell_1$ consecutive 0's followed by (at least $\ell_2$) at most $k_2$ consecutive 1's. Here, $r= \ell_1+k_2$, ${\mathbb P}\big(H_{\ell_1,\ell_2,k_2}^0=0\big)=1$ and $\varepsilon_n:=\sup\big\{x:{\mathbb P}\big(H_{\ell_1,\ell_2,k_2}^n=x\big)>0\big\}=\floor{n/\ell}$. Also, if $0$ occurs after at least $\ell_1$ consecutive $0$'s followed by (at least $\ell_2$) at most $k_2$ consecutive $1$'s then $H_{\ell_1,\ell_2,k_2}^n$ moves $v$ (any) to $v+1$ times. Therefore, $\kappa_0=(1,0,\dotsc,0)_{1\times (\ell_1+k_2+1)}$, $A=[a_{i,j}]_{(\ell_1+k_2+1)\times(\ell_1+k_2+1)}$ with non-zero entries
\begin{itemize}
\item $a_{i,1}=p$ and $a_{i,i+1}=q$ for $1 \le i \le \ell_1$,
\item $a_{\ell_1+1,\ell_1+1}=q$ and $a_{\ell_1+1,\ell_1+2}=p$,
\item $a_{i,2}=q$ for $\ell_1+2 \le i \le \ell_1+\ell_2$ and $a_{i,i+1}=p$ for $\ell_1+2 \le i \le \ell_1+k_2$,
\item $a_{\ell_1+k_2+1,1}=p$
\end{itemize}
and $B=[b_{i,j}]_{(\ell_1+k_2+1)\times(\ell_1+k_2+1)}$ is the matrix of non-zero entries $b_{i,2}=q$ for $\ell_1+\ell_2+1 \le i \le \ell_1+k_2+1$. Hence, using \eqref{four:dpgf}, it can be easily verified that
\begin{equation}
\Phi(t,z)= \frac{1-a(p)z^\ell(t-1)\displaystyle{\sum_{i=1}^{m_2}(pz)^{i-1}}}{1-z-a(p)z^\ell(t-1)\left(1-(pz)^{m_2}\right)}.\label{four:ddpgf1}
\end{equation}
Now, using \eqref{four:ddpgf1}, the following theorem can be easily derived.
\begin{theorem}\label{four:thh1}
The recursive relation in PGF, PMF and moments of $H_{\ell_1,\ell_2,k_2}^n$, for $n \ge \ell+1$, are given by
\begin{itemize}
\item[(i)] $\phi_n(t)=\phi_{n-1}(t)+a(p)(t-1)\left[\phi_{n-\ell}(t)-p^{m_2}\phi_{n-\ell-m_2}(t)\right]-a(p)(t-1)p^{n-\ell}{\bf 1}(\ell+1 \le n\le \ell+m_2-1)$\\
with initial condition $\phi_n(t)=1,$ for $n \le \ell$, where ${\bf 1}(A)$ denotes the indicator function of set $A$.
\item[(ii)] $\begin{aligned}[t]
p_{m,n}&=p_{m,n-1}+a(p)\left[p_{m-1,n-\ell}-p_{m,n-\ell}-p^{m_2}(p_{m-1,n-\ell-m_2}-p_{m,n-\ell-m_2})\right]\\
&~~~-a(p)~p^{n-\ell}\left[{\bf 1}(m=1,~\ell+1\le n \le \ell+ m_2-1)-{\bf 1}(m=0,~\ell+1\le n \le \ell+m_2-1)\right]\\
\end{aligned}$\\
with initial conditions $p_{0,n} =1,$ $p_{m,n}=0,~m>0$ for $n \le \ell$.
\item[(iii)] $\begin{aligned}[t]
\mu_{n,j}&=\mu_{n-1,j}+a(p)\sum_{k=0}^{j-1}{j \choose k}\left[\mu_{n-\ell,k}-p^{m_2}\mu_{n-\ell-m_2,k}\right]-a(p)~p^{n-\ell}{\bf 1}(\ell+1 \le n \le\ell+m_2-1),\\
\end{aligned}$\\
for $j \ge 1$ with initial conditions $\mu_{n,0}=1$ and $\mu_{n,j}=0$ for all $j \ge 1$ and $n \le \ell$.
\end{itemize}
\end{theorem}
\noindent
Next, we obtain an explicit form for PGF and PMF using Theorem \ref{four:thh1}.
\begin{theorem}\label{four:ex2}
Assume the conditions of Theorem \ref{four:thh1} hold, then PGF and PMF of $H_{\ell_1,\ell_2,k_2}^n$ are given by
\begin{itemize}
\item[(i)] $\phi_n(t) = \chi_n(t)-a(p)(t-1)\sum_{i=\ell}^{\ell+m_2-1}p^{i-\ell}\chi_{n-i}(t)$
\item[(ii)] $p_{m,n}={\cal V}_{m,n}-a(p)\sum_{i=\ell}^{\ell+m_2-1}p^{i-\ell}({\cal V}_{m-1,n-i}-{\cal V}_{m,n-i}),$
\end{itemize}
where
$$\chi_n(t)=\sum_{u=0}^{\floor{\frac{n}{\ell}}}\sum_{v=0}^{\floor{\frac{n-u \ell}{\ell+m_2}}}{n-u (\ell-1)-v(\ell+m_2-1) \choose n-u \ell-v(\ell+m_2),u,v}(-1)^v p^{v m_2}(a(p)(t-1))^{u+v}$$
and
$${\cal V}_{m,n}=\sum_{u=0}^{\floor{\frac{n}{\ell}}}\sum_{v=0}^{\floor{\frac{n-u \ell}{\ell+m_2}}}{n-u (\ell-1)-v(\ell+m_2-1) \choose n-u \ell-v(\ell+m_2),u,v}{u+v \choose m}(-1)^{u-m} p^{v m_2} a(p)^{u+v}.$$
\end{theorem}
\noindent
Next, using \eqref{four:wtdpgf}, it can be easily verified that
\begin{equation}
H(t,z)=1+\frac{qt}{1-pt}\sum_{r=1}^{\infty}\left(\frac{a(p)t^\ell(1-(pt)^{m_2})}{1-t+a(p)t^\ell(1-(pt)^{m_2})}\right)^r z^r.\label{four:wt1}
\end{equation}
Hence, using \eqref{four:wt1}, the following theorem can be easily derived.
\begin{theorem}\label{four:Wth2}
The PGF, PMF and moments of $\rho_{r}$, for $r \ge 1$, are given by
\begin{itemize}
\item[(i)] $\begin{aligned}[t]
H_r(t)=\frac{qt}{1-pt}\left(\frac{a(p)t^\ell(1-(pt)^{m_2})}{1-t+a(p)t^\ell(1-(pt)^{m_2})}\right)^r.
\end{aligned}$
\item[(ii)] $g_r(m)=g_r(m-1)+a(p)[g_{r-1}(m-\ell)-g_r(m-\ell)-p^{m_2}(g_{r-1}(m-\ell-m_2)-g_r(m-\ell-m_2))],$~~$r \ge 2$\\
with initial condition $g_0(m)=\delta_{m,0}$ and\\
$g_1(m)=g_1(m-1)-a(p)[g_1(m-\ell)-p^{m_2}g_1(m-\ell-m_2)]+q a(p) p^{m-\ell-1}{\bf 1}(\ell+1 \le m \le \ell+m_2),$\\
for $m \ge \ell r+1$, $g_r(m)=0$ whenever $m \le \ell r$ and $r \ge 1$.
\item[(iii)] $\tilde{\mu}_{r,j}=\displaystyle{\sum_{k=0}^{j}{j \choose k}\left[\tilde{\mu}_{r,k}+a(p)\left(\ell^{j-k}-p^{m_2}(\ell+m_2)^{j-k}\right)(\tilde{\mu}_{r-1,k}-\tilde{\mu}_{r,k})\right]},$ $j \ge 1$ and $r \ge 2$\vspace{0.25cm}\\
with initial condition $\tilde{\mu}_{0,i}=\delta_{i,0}$ and\\
$\tilde{\mu}_{1,j}=\displaystyle{\sum_{k=0}^{j}{j \choose k}\tilde{\mu}_{1,k}\left[1-a(p)\left(\ell^{j-k}-p^{m_2}(\ell+m_2)^{j-k}\right)\right]+q a(p)\sum_{k=\ell+1}^{\ell+m_2}k^j p^{k-\ell-1}}.$
\end{itemize}
\end{theorem}
\noindent
The proofs of Theorems \ref{four:thh1} - \ref{four:Wth2} follow using steps similar to the proofs of Theorems \ref{four:thh2} - \ref{four:Wth1}.

\subsection{Distribution of $H_{\ell_1,k_1,\ell_2,k_2}^n$ and its Waiting Time}\label{four:D3}
Recall that $H_{\ell_1,k_1,\ell_2,k_2}^n$ is the number of occurrences of (at least $\ell_1$) at most $k_1$ consecutive 0's followed by (at least $\ell_2$) at most $k_2$ consecutive 1's. Here, $r= k_1+k_2+1$ and $k_1^+$ is the element after $k_1$ consecutive 0's (if failures occur) in $\{0,1,\dotsc,k_1,k_1^+=k_1+1,k_1+2,\dotsc,k_1+k_2+1\}$. It is easy to see that ${\mathbb P}\big(H_{\ell_1,k_1,\ell_2,k_2}^0=0\big)=1$ and $\varepsilon_n:=\sup\big\{x:{\mathbb P}\big(H_{\ell_1,k_1,\ell_2,k_2}^n=x\big)>0\big\}=\floor{n/\ell}$. Also, if $0$ occurs after (at least $\ell_1$) at most $k_1$ consecutive $0$'s followed by (at least $\ell_2$) at most $k_2$ consecutive $1$'s then $H_{\ell_1,k_1,\ell_2,k_2}^n$ moves $v$ (any) to $v+1$ times. Therefore, $\kappa_0=(1,0,\dotsc,0)_{1\times (k_1+k_2+2)}$, $A=[a_{i,j}]_{(k_1+k_2+2)\times(k_1+k_2+2)}$ with non-zero entries
\begin{itemize}
\item $a_{i,1}=p$ and $a_{i,i+1}=q$ for $1 \le i \le \ell_1$,
\item $a_{i,k_1+3}=p$ and $a_{i,i+1}=q$ for $\ell_1+1 \le i \le k_1+1$,
\item $a_{k_1+2,1}=p$ and $a_{k_1+2,k_1+2}=q$,
\item $a_{i,2}=q$ for $k_1+3 \le i \le k_1+\ell_2+1$ and $a_{i,i+1}=p$ for $k_1+3 \le i \le k_1+k_2+1$,
\item $a_{k_1+k_2+2,1}=p$
\end{itemize}
and $B=[b_{i,j}]_{(k_1+k_2+2)\times(k_1+k_2+2)}$ is the matrix of non-zero entries $b_{i,2}=q$ for $k_1+\ell_2+2 \le i \le k_1+k_2+2$. Hence, using \eqref{four:dpgf}, it can be easily verified that
\begin{equation}
\Phi(t,z)=\sum_{n=0}^{\infty}\phi_n(t)z^n= \frac{1-a(p)z^\ell(t-1)\left(1-(qz)^{m_1}\right)\displaystyle{\sum_{i=1}^{m_2}(pz)^{i-1}}}{1-z-a(p)z^\ell(t-1)\left(1-(qz)^{m_1}\right)\left(1-(pz)^{m_2}\right)}.\label{four:ddpgf3}
\end{equation}
Now, using \eqref{four:ddpgf3}, the following theorem can be easily derived.
\begin{theorem}\label{four:th1}
The recursive relations in PGF, PMF and moments of $H_{\ell_1,k_1,\ell_2,k_2}^n$, for $n \ge \ell+1$, are given by
\begin{itemize}
\item[(i)] $\begin{aligned}[t]
\phi_n(t)&=\phi_{n-1}(t)+a(p)(t-1)\left[\phi_{n-\ell}(t)-q^{m_1}\phi_{n-\ell-m_1}(t)-p^{m_2}\phi_{n-\ell-m_2}(t)+q^{m_1}p^{m_2}\phi_{n-\ell-m_1-m_2}(t)\right]\\
&~~~-a(p)(t-1)p^{n-\ell}\left({\bf 1}(\ell+1 \le n\le \ell+m_2-1)-\left(\frac{q}{p}\right)^{m_1}{\bf 1}(\ell+m_1\le n \le\ell+m_1+m_2-1)\right)
\end{aligned}$\\
\\
with initial condition $\phi_n(t)=1,$ for $n \le \ell$.
\item[(ii)] $\begin{aligned}[t]
p_{m,n}&=p_{m,n-1}-a(p)~p^{n-\ell}\Bigr[{\bf 1}(m=1,~\ell+1\le n \le \ell+ m_2-1)-{\bf 1}(m=0,~\ell+1\le n \le \ell+m_2-1)\\
&~~~-\left.\left(q/p\right)^{m_1}\Big({\bf 1}(m\hspace{-0.02cm}=\hspace{-0.02cm}1,~\ell\hspace{-0.02cm}+\hspace{-0.02cm}m_1 \hspace{-0.02cm}\le\hspace{-0.02cm} n\hspace{-0.02cm} \le\hspace{-0.02cm}\ell\hspace{-0.02cm}+\hspace{-0.02cm}m_1\hspace{-0.02cm}+\hspace{-0.02cm}m_2\hspace{-0.02cm}-\hspace{-0.02cm}1)\hspace{-0.07cm}-\hspace{-0.07cm}{\bf 1}(m\hspace{-0.02cm}=\hspace{-0.02cm}0,~\ell\hspace{-0.02cm}+\hspace{-0.02cm}m_1 \hspace{-0.02cm}\le\hspace{-0.02cm}n\hspace{-0.02cm} \le\hspace{-0.02cm} \ell\hspace{-0.02cm}+\hspace{-0.02cm}m_1\hspace{-0.02cm}+\hspace{-0.02cm}m_2\hspace{-0.02cm}-\hspace{-0.02cm}1)\Bigr)\right]\\
&~~~+a(p)\left[p_{m-1,n-\ell}-p_{m,n-\ell}-q^{m_1}(p_{m-1,n-\ell-m_1}-p_{m,n-\ell-m_1})-p^{m_2}(p_{m-1,n-\ell-m_2}-p_{m,n-\ell-m_2})\right.\\
&~~~\left.+q^{m_1}p^{m_2}(p_{m-1,n-\ell-m_1-m_2}-p_{m,n-\ell-m_1-m_2})\right]
\end{aligned}$\\
\\
with initial conditions $p_{0,n} =1$ and $p_{m,n}=0,~m>0$ for $n \le \ell$.
\item[(iii)] $\begin{aligned}[t]
\mu_{n,j}&=\mu_{n-1,j}+a(p)\sum_{k=0}^{j-1}{j \choose k}\left[\mu_{n-\ell,k}-q^{m_1}\mu_{n-\ell-m_1,k}-p^{m_2}\mu_{n-\ell-m_2,k}+q^{m_1}p^{m_2}\mu_{n-\ell-m_1-m_2,k}\right]\\
&~~~-a(p)~p^{n-\ell}\Big[{\bf 1}(\ell+1 \le n \le\ell+m_2-1)-\left(q/p\right)^{m_1}{\bf 1}(\ell+m_1 \le n \le\ell+m_1+m_2-1)\Big],~~j \ge 1
\end{aligned}$\\
\\
with initial conditions $\mu_{n,0}=1$ and $\mu_{n,j}=0$ for all $j \ge 1$ and $n \le \ell$.
\end{itemize}
\end{theorem}
\noindent
Next, we obtain an explicit form for PGF and PMF using Theorem \ref{four:th1}.
\begin{theorem}\label{four:ex3}
Assume the conditions of Theorem \ref{four:th1} hold, then PGF and PMF of $H_{\ell_1,k_1,\ell_2,k_2}^n$ are given by
\begin{itemize}
\item[(i)] $\begin{aligned}[t]
\phi_n(t) = \varphi_n(t)-a(p)(t-1)\left[\sum_{i=\ell}^{\ell+m_2-1}p^{i-\ell}\varphi_{n-i}(t)-\left(\frac{q}{p}\right)^{m_1}\sum_{i=\ell+m_1}^{\ell+m_1+m_2-1}p^{i-\ell}\varphi_{n-i}(t)\right]
\end{aligned}$
\item[(ii)] $\begin{aligned}[t]
p_{m,n}=\kappa_{m,n}-a(p)\left[\sum_{i=\ell}^{\ell+m_2-1}p^{i-\ell}(\kappa_{m-1,n-i}-\kappa_{m,n-i})-\left(\frac{q}{p}\right)^{m_1}\sum_{i=\ell+m_1}^{\ell+m_1+m_2-1}p^{i-\ell}(\kappa_{m-1,n-i}-\kappa_{m,n-i})\right],
\end{aligned}$
\end{itemize}
where
\begin{align*}
\varphi_n(t)\hspace{-0.1cm}&=\hspace{-0.15cm}\sum_{u=0}^{\floor{\frac{n}{\ell}}}\hspace{-0.1cm}\sum_{w=0}^{\floor{\frac{n-u \ell}{\ell+m_1}}}\hspace{-0.1cm}\sum_{r=0}^{\floor{\frac{f(n,u,w,0,0)}{\ell+m_2}}}\hspace{-0.1cm}\sum_{v=0}^{\floor{\frac{f(n,u,w,r,0)}{\ell+m_1+m_2}}}\hspace{-0.65cm}(-1)^{w+r}{f(n,u,w,r,v)\hspace{-0.07cm}+\hspace{-0.07cm}u\hspace{-0.07cm}+\hspace{-0.07cm}v\hspace{-0.07cm}+\hspace{-0.07cm}r\hspace{-0.07cm}+\hspace{-0.07cm}w \choose f(n,u,w,r,v),u,w,r,v} q^{(v\hspace{-0.01cm}+\hspace{-0.01cm}w)m_1}p^{(v\hspace{-0.01cm}+\hspace{-0.01cm}r)m_2}(a(p)(t\hspace{-0.07cm}-\hspace{-0.07cm}1))^{u\hspace{-0.01cm}+\hspace{-0.01cm}w\hspace{-0.01cm}+\hspace{-0.01cm}r\hspace{-0.01cm}+\hspace{-0.01cm}v}
\end{align*}
\begin{align*}
\kappa_{m,n}&=\sum_{u=0}^{\floor{\frac{n}{\ell}}} \sum_{w=0}^{\floor{\frac{n-u \ell}{\ell+m_1}}}\sum_{r=0}^{\floor{\frac{f(n,u,w,0,0)}{\ell+m_2}}} \sum_{v=0}^{\floor{\frac{f(n,u,w,r,0)}{\ell+m_1+m_2}}} (-1)^{u+v-m}{f(n,u,w,r,v)\hspace{-0.07cm}+\hspace{-0.07cm}u\hspace{-0.07cm}+\hspace{-0.07cm}v\hspace{-0.07cm}+\hspace{-0.07cm}r\hspace{-0.07cm}+\hspace{-0.07cm}w \choose f(n,u,w,r,v),u,w,r,v}{u+w+r+v \choose m} q^{(v+w)m_1}\\
&~~~~~~~~~~~~~~~~~~~~~~~~~~~~~~~~~~~~~~~~~~~~~~~~~~~~~~~~~~~~~~~~~~~~~~~~~~~~~~~~~~~~~~~~~~~~~~~~~~~~~p^{(v+r)m_2} a(p)^{u+w+r+v}
\end{align*}
and $f(n,u,w,r,v)=n-u\ell-w(\ell+m_1)-r(\ell+m_2)-v(\ell+m_1+m_2)$.
\end{theorem}
\noindent
Next, using \eqref{four:wtdpgf}, it can be easily verified that
\begin{equation}
H(t,z)=1+\frac{qt}{1-pt}\sum_{r=1}^{\infty}\left(\frac{a(p)t^\ell(1-(qt)^{m_1})(1-(pt)^{m_2})}{1-t+a(p)t^\ell(1-(qt)^{m_1})(1-(pt)^{m_2})}\right)^r z^r.\label{four:wt3}
\end{equation}
Hence, using \eqref{four:wt3}, the following theorem can be easily derived.
\begin{theorem}\label{four:wth3}
The PGF, PMF and moments of $\rho_{r}$, for $r \ge 1$, are given by
\begin{itemize}
\item[(i)] $\begin{aligned}[t]
H_r(t)=\frac{qt}{1-pt}\left(\frac{a(p)t^\ell(1-(qt)^{m_1})(1-(pt)^{m_2})}{1-t+a(p)t^\ell(1-(qt)^{m_1})(1-(pt)^{m_2})}\right)^r.
\end{aligned}$
\item[(ii)] $\begin{aligned}[t]
g_r(m)&\hspace{-0.07cm}=\hspace{-0.07cm}g_r(m-1)+a(p)[g_{r-1}(m-\ell)-g_r(m-\ell)-q^{m_1}(g_{r-1}(m-\ell-m_1)-g_r(m-\ell-m_1))\\
&~~~-p^{m_2}(g_{r-1}(m\hspace{-0.06cm}-\hspace{-0.06cm}\ell\hspace{-0.06cm}-\hspace{-0.06cm}m_2)\hspace{-0.06cm}-\hspace{-0.06cm}g_r(m\hspace{-0.06cm}-\hspace{-0.06cm}\ell\hspace{-0.06cm}-\hspace{-0.06cm}m_2))\hspace{-0.06cm}+\hspace{-0.06cm}q^{m_1}p^{m_2}(g_{r-1}(m\hspace{-0.06cm}-\hspace{-0.06cm}\ell\hspace{-0.06cm}-\hspace{-0.06cm}m_1\hspace{-0.06cm}-\hspace{-0.06cm}m_2)\hspace{-0.06cm}-\hspace{-0.06cm}g_r(m\hspace{-0.06cm}-\hspace{-0.06cm}\ell\hspace{-0.06cm}-\hspace{-0.06cm}m_1\hspace{-0.06cm}-\hspace{-0.06cm}m_2))],
\end{aligned}$\\
\\
for $r \ge 2$ with initial condition $g_0(m)=\delta_{m,0}$ and\vspace{-0.43cm}
\begin{align*}
g_1(m)&=g_1(m-1)\hspace{-0.02cm}+\hspace{-0.02cm}q a(p) p^{m-\ell-1}\left({\bf 1}(\ell\hspace{-0.02cm}+\hspace{-0.02cm}1 \hspace{-0.02cm}\le\hspace{-0.02cm} m \hspace{-0.02cm}\le\hspace{-0.02cm} \ell\hspace{-0.02cm}+\hspace{-0.02cm}m_2)\hspace{-0.02cm}-\hspace{-0.02cm}\left(\frac{q}{p}\right)^{m_1}{\bf 1}(\ell\hspace{-0.02cm}+\hspace{-0.02cm}m_1\hspace{-0.02cm}+\hspace{-0.02cm}1\hspace{-0.02cm}\le\hspace{-0.02cm} m\hspace{-0.02cm}\le\hspace{-0.02cm} \ell\hspace{-0.02cm}+\hspace{-0.02cm}m_1\hspace{-0.02cm}+\hspace{-0.02cm}m_2)\right)\\
&~~~-a(p)[g_1(m-\ell)-q^{m_1}g_1(m-\ell-m_1)-p^{m_2}g_1(m-\ell-m_2)+q^{m_1}p^{m_2}g_1(m-\ell-m_1-m_2)],
\end{align*}\vskip -5ex
for $m \ge \ell r+1$, $g_r(m)=0$ whenever $m \le \ell r$ and $r \ge 1$.
\item[(iii)] $\begin{aligned}[t]
\tilde{\mu}_{r,j}&=\displaystyle{\sum_{k=0}^{j}{j \choose k}\left[\tilde{\mu}_{r,k}+a(p)\left(\ell^{j-k}-q^{m_1}(\ell+m_1)^{j-k}-p^{m_2}(\ell+m_2)^{j-k}\right.\right.}\\
&~~~+\left.\left.q^{m_1}p^{m_2}(\ell+m_1+m_2)^{j-k}\right)\right](\tilde{\mu}_{r-1,k}-\tilde{\mu}_{r,k}),~~~ j \ge 1,~{\rm and}~ r \ge 2
\end{aligned}$ \vspace{0.25cm}\\
with initial condition $\tilde{\mu}_{0,i}=\delta_{i,0}$ and\\
$\begin{aligned}[t]
\tilde{\mu}_{1,j}&=\sum_{k=0}^{j}{j \choose k}\tilde{\mu}_{1,k}\left[1-a(p)\left(\ell^{j-k}-q^{m_1}(\ell+m_1)^{j-k}-p^{m_2}(\ell+m_2)^{j-k}+q^{m_1}p^{m_2}(\ell+m_1+m_2)^{j-k}\right)\right]\\
&~~~+q a(p)\left[\sum_{k=\ell+1}^{\ell+m_2}k^j p^{k-\ell-1} - \left(\frac{q}{p}\right)^{m_1}\sum_{k=\ell+m_1+1}^{\ell+m_1+m_2}k^j p^{k-\ell-1}\right].
\end{aligned}$
\end{itemize}
\end{theorem}
\noindent
The proofs of Theorems \ref{four:th1} - \ref{four:wth3} follow using steps similar to the proofs of Theorems \ref{four:thh2} - \ref{four:Wth1}.

\subsection{Generalization of $H_{\ell_1,k_1,\ell_2}$, $H_{\ell_1,\ell_2,k_2}$ and $H_{\ell_1,k_1,\ell_2,k_2}$ under Markov Dependence}\label{four:MDR}
We now demonstrate that the results for $H_{\ell_1,k_1,\ell_2}$, $H_{\ell_1,\ell_2,k_2}$ and $H_{\ell_1,k_1,\ell_2,k_2}$ can be effortlessly generalized to Markov dependent setup. Let $\xi_1,\xi_2,\dotsc$ be time homogeneous two-state Markov chain with transition probability matrix 
$$
P=\begin{pmatrix}
p_{00}&p_{01}\\
p_{10}&p_{11}
\end{pmatrix},
$$
i.e., $p_{rs}={\mathbb P}(\xi_n=s|\xi_{n-1}=r)$ for $n\ge2$ and $r,s \in \{0,1\}$ and initial distribution $p_i={\mathbb P}(\xi_1=i)$, $i=0,1$.\\
Now, we derive the double generating function for $H_{\ell_1,k_1,\ell_2}$, $H_{\ell_1,\ell_2,k_2}$ and $H_{\ell_1,k_1,\ell_2,k_2}$, and their waiting time using Markov chain approach as discussed in Section \ref{four:DR}. The expressions for double generating functions follow directly from \eqref{four:dpgf} and \eqref{four:wtdpgf} with $\kappa_0$, $A$ and $B$ as defined below. 

\begin{itemize}
\item[{\bf (T1)}] Here, $\kappa_0=(p_0,p_1,0,\dotsc,0)_{1\times (k_1+\ell_2+2)}$, $A=[a_{i,j}]_{(k_1+\ell_2+2)\times(k_1+\ell_2+2)}$ with non-zero entries
\begin{itemize}
\item[$\bullet$] $a_{11}=p_{11}$, $a_{12}=p_{10}$, $a_{i,1}=p_{01}$ and $a_{i,i+1}=p_{00}$ for $2 \le i \le \ell_1$,
\item[$\bullet$] $a_{i,k_1+3}=p_{01}$, $a_{i,i+1}=p_{00}$ for $\ell_1+1 \le i \le k_1+1$, $a_{k_1+2,1}=p_{01}$ and $a_{k_1+2,k_1+2}=p_{00}$,
\item[$\bullet$] $a_{i,2}\hspace{-0.04cm}=\hspace{-0.04cm}p_{10}$ for $k_1+3 \hspace{-0.04cm}\le\hspace{-0.04cm} i \hspace{-0.04cm}\le\hspace{-0.04cm} k_1+\ell_2+2$, $a_{i,i+1}\hspace{-0.04cm}=\hspace{-0.04cm}p_{11}$ for $k_1+3 \hspace{-0.04cm}\le\hspace{-0.04cm} i \hspace{-0.04cm}\le\hspace{-0.04cm} k_1+\ell_2$ and $a_{k_1+\ell_2+2,k_1+\ell_2+2}\hspace{-0.04cm}=\hspace{-0.04cm}p_{11}\hspace{-0.02cm},$
\end{itemize}
and $B=[b_{i,j}]_{(k_1+\ell_2+2)\times(k_1+\ell_2+2)}$ is the matrix of non-zero entry $b_{k_1+\ell_2+1,k_1+\ell_2+2}=p_{11}$. Hence, using \eqref{four:dpgf} and \eqref{four:wtdpgf}, it can be verified that 
\begin{equation}
\Phi(t,z)=\frac{p_0[1+(p_{10}-p_{00})z]+p_1[1+(p_{01}-p_{11})z+A(z)(t-1)(1-(p_{00}z)^{m_1})]}{1-(p_{00}+p_{11})z-(p_{01}p_{10}-p_{00}p_{11})z^2-A(z)(t-1)(1-(p_{00}z)^{m_1})}\label{four:Mddpgf2}
\end{equation}
and
\begin{equation}
H(t,z)=1+\frac{p_0 p_{10}t+p_1(1-p_{11}t)}{p_{10}t}\sum_{r=1}^{\infty}\left(\frac{A(t)(1-(p_{00}t)^{m_1})}{1-(p_{00}+p_{11})t-(p_{01}p_{10}-p_{00}p_{11})t^2+A(t)(1-(p_{00}t)^{m_1})}\right)^r \hspace{-0.14cm}z^r\hspace{-0.1cm},\label{four:Mwt2}
\end{equation}
where $A(t)=(p_{01}t)(p_{10}t)(p_{00}t)^{\ell_1-1}(p_{11}t)^{\ell_2-1}$.

\item[{\bf (T2)}] Here, $\kappa_0=(p_0,p_1,0,\dotsc,0)_{1\times (\ell_1+k_2+1)}$, $A=[a_{i,j}]_{(\ell_1+k_2+1)\times(\ell_1+k_2+1)}$ with non-zero entries
\begin{itemize}
\item[$\bullet$] $a_{11}=p_{11}$, $a_{12}=p_{10}$, $a_{i,1}=p_{10}$, $a_{i,i+1}=p_{00}$ for $2 \le i \le \ell_1$, $a_{\ell_1+1,\ell_1+1}=p_{00}$ and $a_{\ell_1+1,\ell_1+2}=p_{01}$,
\item[$\bullet$] $a_{i,2}=p_{10}$ for $\ell_1+2 \le i \le \ell_1+\ell_2$, $a_{i,i+1}=p_{11}$ for $\ell_1+2 \le i \le \ell_1+k_2$ and $a_{\ell_1+k_2+1,1}=p_{11},$
\end{itemize}
and $B=[b_{i,j}]_{(\ell_1+k_2+1)\times(\ell_1+k_2+1)}$ is the matrix of non-zero entries $b_{i,2}=p_{10}$ for $\ell_1+\ell_2+1 \le i \le \ell_1+k_2+1$. Hence, using \eqref{four:dpgf} and \eqref{four:wtdpgf}, it can be easily verified that
\begin{equation}
\Phi(t,z)= \frac{p_0\Big[1+(p_{10}-p_{00})z-A(z)(t-1)\displaystyle{\sum_{i=1}^{m_2}(pz)^{i-1}}\Bigr]+p_1[1+(p_{01}-p_{11})z]}{1-(p_{00}+_{11})z-(p_{01}p_{10}-p_{00}p_{11})z^2-A(z)(t-1)\left(1-(p_{11}z)^{m_2}\right)}\label{four:Mddpgf1}
\end{equation}
and
\begin{equation}
H(t,z)=1+\frac{p_0 p_{10}t+p_1(1-p_{11}t)}{1-p_{11}t}\sum_{r=1}^{\infty}\left(\frac{A(t)(1-(p_{11}t)^{m_2})}{1-(p_{00}+p_{11})t-(p_{01}p_{10}-p_{00}p_{11})t^2+A(t)(1-(p_{11}t)^{m_2})}\right)^r\hspace{-0.15cm} z^r\hspace{-0.1cm}.\label{four:MMMwt1}
\end{equation}

\item[{\bf (T3)}] Here, $\kappa_0=(p_0,p_1,0,\dotsc,0)_{1\times (k_1+k_2+2)}$, $A=[a_{i,j}]_{(k_1+k_2+2)\times(k_1+k_2+2)}$ with non-zero entries
\begin{itemize}
\item[$\bullet$] $a_{11}=p_{11}$, $a_{12}=p_{10}$, $a_{i,1}=p_{01}$ and $a_{i,i+1}=p_{00}$ for $2 \le i \le \ell_1$,
\item[$\bullet$] $a_{i,k_1+3}=p_{01}$, $a_{i,i+1}=p_{00}$ for $\ell_1+1 \le i \le k_1+1$, $a_{k_1+2,1}=p_{01}$ and $a_{k_1+2,k_1+2}=p_{00}$,
\item[$\bullet$] $a_{i,2}=p_{10}$ for $k_1+3 \le i \le k_1+\ell_2+1$, $a_{i,i+1}=p_{11}$ for $k_1+3 \le i \le k_1+k_2+1$ and $a_{k_1+k_2+2,1}=p_{11},$
\end{itemize}
and $B=[b_{i,j}]_{(k_1+k_2+2)\times(k_1+k_2+2)}$ is the matrix of non-zero entries $b_{i,2}=p_{10}$ for $k_1+\ell_2+2 \le i \le k_1+k_2+2$. Hence, using \eqref{four:dpgf} and \eqref{four:wtdpgf}, it can be easily verified that
\begin{equation}
\Phi(t,z)= \frac{p_0\Big[1+(p_{10}-p_{00})z-A(z)(t-1)(1-(p_{00}z)^{m_1})\displaystyle{\sum_{i=1}^{m_2}(pz)^{i-1}}\Bigr]+p_1[1+(p_{01}-p_{11})z]}{1-(p_{00}+p_{11})z-(p_{01}p_{10}-p_{00}p_{11})z^2-A(z)(t-1)\left(1-(p_{00}z)^{m_1}\right)\left(1-(p_{11}z)^{m_2}\right)}\label{four:Mddpgf3}
\end{equation}
and
\begin{equation}
H(t,z)\hspace{-0.1cm}=\hspace{-0.1cm}1+\frac{p_0 p_{10}t\hspace{-0.05cm}+\hspace{-0.05cm}p_1(1\hspace{-0.05cm}-\hspace{-0.05cm}p_{11}t)}{1-p_{11}t}\hspace{-0.1cm}\sum_{r=1}^{\infty}\hspace{-0.1cm}\left(\hspace{-0.1cm}\frac{A(t)(1-(p_{00}t)^{m_1})(1-(p_{11}t)^{m_2})z}{1\hspace{-0.08cm}-\hspace{-0.08cm}(p_{00}\hspace{-0.08cm}+\hspace{-0.08cm}p_{11})t\hspace{-0.08cm}-\hspace{-0.08cm}(p_{01}p_{10}\hspace{-0.08cm}-\hspace{-0.08cm}p_{00}p_{11})t^2\hspace{-0.08cm}+\hspace{-0.08cm}A(t)(1\hspace{-0.08cm}-\hspace{-0.08cm}(p_{00}t)^{m_1})(1\hspace{-0.08cm}-\hspace{-0.08cm}(p_{11}t)^{m_2})}\hspace{-0.1cm}\right)^r \hspace{-0.2cm}.\label{four:MMMwt3}
\end{equation}
\end{itemize}

\section{Applications}\label{APP}
In this section, we discuss the relevance of the results derived through Fibonacci words and quality control. Also, for more applications, we refer the reader to Balakrishnan and Koutras \cite{BK}, Fu \cite{FU},  Balakrishnan {\em et al.} \cite{BKM1, BKM2}, Moore \cite{moore} and references therein.

\subsection{Fibonacci Words}
Fibonacci words are particular sequences of binary numbers $0$ and $1$ (or two alphabets) and it is used to model physical systems with the aperiodic order such as quasi-crystals. Also, Fibonacci word have been studied widely in the field of combinatorics on words. Fibonacci words are formed in a similar way as Fibonacci numbers (repeated addition) and, in this process, $n$-th Fibonacci word depends on $(n-1)$-th and $(n-2)$-th Fibonacci words of $0$'s and $1$'s. The construction can be explained as follows:
$${\cal C}_0=0\quad\quad {\rm and}\quad \quad {\cal C}_1=01$$
\noindent
then $n$-th Fibonacci word is given by
$${\cal C}_n={\cal C}_{n-1} {\cal C}_{n-2}.$$
For example, 10-th element of Fibonacci words is given by\\
${\cal C}_{10}$ = 0 1 0 0 1 0 1 0 0 1 0 0 1 0 1 0 0 1 0 1 0 0 1 0 0 1 0 1 0 0 1 0 0 1 0 1 0 0 1 0 1 0 0 1 0 0 1 0 1 0 0 1 0 1 0 0 1 0 0 1 0 1 0 0 1 0 0 1 0 1 0 0 1 0 1 0 0 1 0 0 1 0 1 0 0 1 0 0 1 0 1 0 0 1 0 1 0 0 1 0 0 1 0 1 0 0 1 0 1 0 0 1 0 0 1 0 1 0 0 1 0 0 1 0 1 0 0 1 0 1 0 0 1 0 0 1 0 1 0 0 1 0 1 0\\
and the random variable representation is given by\\
$(1-\zeta_{1})\zeta_{2}(1-\zeta_{3})(1-\zeta_{4})\zeta_{5}(1-\zeta_{6})\zeta_{7}(1-\zeta_{8})(1-\zeta_{9})\zeta_{10}(1-\zeta_{11})(1-\zeta_{12})\zeta_{13}(1-\zeta_{14})\zeta_{15}(1-\zeta_{16})(1-\zeta_{17})\zeta_{18}(1-\zeta_{19})\zeta_{20}(1-\zeta_{21})(1-\zeta_{22})\zeta_{23}(1-\zeta_{24})(1-\zeta_{25})\zeta_{26}(1-\zeta_{27})\zeta_{28}(1-\zeta_{29})(1-\zeta_{30})\dotsc$.\\
Also, the sub-words ``11'' and ``000'' never occur in Fibonacci words and last two digits are ``01'' and ``10'', alternately. For more details on Fibonacci words, we refer the reader to Berstel \cite{Ber}. Now, observe that Fibonacci words can be seen as a pattern of either exactly one 1 followed by (at least one) at most two consecutive 0's or (at least one) at most two consecutive 0's followed by exactly one 1 and hence the distribution of patterns adopted the distribution of either $H_{1,1,1,2}^n$ or $H_{1,2,1,1}^n$ respectively, for $n$-th Fibonacci word. For large values of $n$, the probabilities and moments of the distribution of these patterns can be calculated from the distribution of either $H_{1,1,1,2}^n$ or $H_{1,2,1,1}^n$. Next, we compute some probabilities and mean for $H_{1,2,1,1}^n$ and its waiting time for various values of $p$ and $n=60$.\\
\begin{table}[h!]
  \centering
  \caption{Distribution and moments of $H_{1,2,1,1}^{60}$.}
  \label{four:ta4}
  \begin{tabular}{lllllllll}
    \toprule
    $n$ & $m$ & $p=0.35$ & $p=0.36$ & $p=0.37$ & $p=0.38$ & $p=0.39$ & $p=0.40$\\
    \midrule
    \multirow{6}{*}{60}& 0  & 0.0081259 & 0.0073285 & 0.0066661 & 0.0061179 & 0.0056670 & 0.0052998\\
     &                    1 & 0.0363192 & 0.0335666 & 0.0312188 & 0.0292301 & 0.0275615 & 0.0261798\\
     &                    2 & 0.0844787 & 0.0798366 & 0.0757692 & 0.0722423 & 0.0692234 & 0.0666826\\
     &                    3 & 0.1353360 & 0.1305530 & 0.1262260 & 0.1223700 & 0.1189930 & 0.1160990\\
     &                    4 & 0.1669740 & 0.1641700 & 0.1614830 & 0.1589750 & 0.1566960 & 0.1546850\\
     &                    5 & 0.1683560 & 0.1684990 & 0.1684180 & 0.1681850 & 0.1678630 & 0.1675060\\
       \midrule
     ${\mathbb E}\left(H_{1,2,1,1}^{60}\right)$ & -& 5.07803 & 5.17016 & 5.25346 & 5.32777 & 5.39297 & 5.44896\\
     \bottomrule
  \end{tabular}
\end{table}
\begin{table}[h!]
  \centering
  \caption{Distribution and moments of waiting time for $H_{1,2,1,1}^{60}$.}
   \label{four:ta5}
  \begin{tabular}{lllllllll}
    \toprule
    $r$ & $m$ & $p=0.45$ & $p=0.46$ & $p=0.47$ & $p=0.48$ & $p=0.49$ & $p=0.50$\\
    \midrule
    \multirow{8}{*}{1}&  3 & 0.1361250 & 0.1341360 & 0.1320230 & 0.1297920 & 0.1274490 & 0.1250000 \\
     &                   4 & 0.1361250 & 0.1341360 & 0.1320230 & 0.1297920 & 0.1274490 & 0.1250000\\
     &                   5 & 0.0612563 & 0.0617026 & 0.0620508 & 0.0623002 & 0.0624500 & 0.0625000\\
     &                   6 & 0.0427262 & 0.0437101 & 0.0446207 & 0.0454542 & 0.0462068 & 0.0468750\\
     &                   7 & 0.0529177 & 0.0534260 & 0.0538587 & 0.0542141 & 0.0544908 & 0.0546875\\
     &                   8 & 0.0547707 & 0.0548654 & 0.0549045 & 0.0548879 & 0.0548157 & 0.0546875\\
     &                   9 & 0.0464322 & 0.0465889 & 0.0467123 & 0.0468019 & 0.0468565 & 0.0468750\\
     &                  10 & 0.0399053 & 0.0401752 & 0.0404228 & 0.0406466 & 0.0408449 & 0.0410156\\
        \midrule
     ${\mathbb E}\left(\rho_{1}\right)$ & -& 2.17153 & 2.31385 & 2.45255 & 2.58869 & 2.72324 & 2.85714\\
     \bottomrule
  \end{tabular}
\end{table}

\noindent
Observe that the upper range of $m$ is $\floor{n/\ell}=\floor{60/2}=30$, while we obtain the probabilities up to $m=5$ and others can be computed in a similar way. Also, for waiting time distribution, it is known that $m \ge \ell r+1=3$. So, we obtain probabilities by taking $m$ up to $10$ in  Table \ref{four:ta5}. Moment for $H_{1,2,1,1}^{60}$ and $\rho_1$ are obtained in Table \ref{four:ta4} and Table \ref{four:ta5}, respectively.

\subsection{Quality Control}
The quality control is a statistical method which monitors the quality of products and services, and is discussed in \cite{BK,GREEN}. It is also described as statistical process control which uses graphical displays (control charts) to determine a process either to be continued or to be adjusted to achieve the desired quality. Here, we consider the start-up demonstration test with Markov dependence to reject the quality of products.\\
A start-up demonstration test can be considered as a scenario in which a customer is interested to buy certain equipments such as, water pumps, garden tillers, car batteries and power generators among many others. The customer can be accepted/rejected the equipment under some predefined conditions. There are several start-up conditions in which the customer can reject the equipment, for example, (i) $m$ out of $n$ consecutive failures (ii) exactly $k_1$ consecutive successful trials are followed by at least $k_2$ consecutive unsuccessful trials. To fit this in our setting, we proposed a rule that an equipment is rejected if the individual start-ups are Markov dependent with (ii). Therefore, the distribution of our interest becomes $H_{k_1,k_1,k_2}^n$ ({\bf (T1)} type) by changing the role of successes and failures and its waiting time is simply the waiting time at which the customer reject an equipment. For more details, we refer the reader to Balakrishnan {\em et al.} \cite{BKM1, BKM2}.\\
%There is a critical zone on control chart outside of which the process is out of control and this can be defined in various ways according to the problem. In particular, the process can be declared out of control in setting on critical zones such as, $m$ out of $n$ consecutive points fall in critical zone, $m$ points inside and $r$ points outside ($r\ge m$) on the control chart. To fit this in our setting, let the process to be declared out of control if exactly $k_1$ consecutive points are within the control limits and at least $k_2$ consecutive points are outside of the control limits ($k_2\ge k_1$) on the control chart. Therefore, the distribution of our interest becomes $H_{k_1,k_1,k_2}^n$ ({\bf (T1)} type) and its waiting time is simply the waiting time at which the process is declared to be out of control.\\
\noindent
Next, let $\rho_r$ be the $r$-th waiting time for $H_{k_1,k_1,k_2}^n$ under Markov dependent trials. Now, we compute the probabilities for waiting time distribution using the results \eqref{four:Mwt2} with $r=2$, $k_1=2$, $k_2=5$, $p_0=1$ and $p_1=0$, and various values of $p_{00}=1-p_{01}$ and $p_{10}=1-p_{11}$ in Table \ref{three:table2}.

\newpage
\begin{table}[ht!]
  \centering
 \caption{Probabilities for waiting time distribution under Markov dependent trials.}
   \label{three:table2}
 \begin{tabular}{llclllllll}
    \toprule
    \multirow{2}{*}{$r$} & \multirow{2}{*}{$(k_1,k_2)$} &\multirow{2}{*}{$m$} & $p_{00}=0.10$ & $p_{00}=0.30$ & $p_{00}=0.50$ & $p_{00}=0.70$ & $p_{00}=0.90$\\
        &             &    & $p_{11}=0.40$ & $p_{11}=0.60$ & $p_{11}=0.50$ & $p_{11}=0.90$ & $p_{11}=0.10$\\
    \midrule
    \multirow{7}{*}{2} & \multirow{7}{*}{(2,5)} & 14 & $1.9\times 10^{-6}$ & 0.00011851 & 0.00006103 & 0.00018983 & $6.6\times 10^{-11}$\\
    &  & 15 & $1.5\times 10^{-6}$ & 0.00014222 & 0.00006103 & 0.00034170 & $1.3\times 10^{-11}$\\ 
    &  & 16 & $2.9\times 10^{-6}$ & 0.00019436 & 0.00007629 & 0.00047269 & $1.4\times 10^{-11}$\\ 
    &  & 17 & $3.1\times 10^{-6}$ & 0.00024177 & 0.00009155 & 0.00059229 & $1.4\times 10^{-11}$\\ 
    &  & 18 & $4.1\times 10^{-6}$ & 0.00028984 & 0.00010681 & 0.00070573 & $1.5\times 10^{-11}$\\ 
    &  & 19 & $4.7\times 10^{-6}$ & 0.00033782 & 0.00012207 & 0.00081590 & $1.6\times 10^{-11}$\\ 
    &  & 20 & $5.5\times 10^{-6}$ & 0.00038581 & 0.00013732 & 0.00092434 & $1.6\times 10^{-11}$\\ 
   \midrule
\multicolumn{2}{c}{${\mathbb E}\left(\rho_2\right)$}  & -& 2411.270 & 288.696 & 512 & 193.544 & $2.5\times 10^6$\\
   \midrule
\multicolumn{2}{c}{$\mathrm{Var}\left(\rho_2\right)$} & -& $2.8\times 10^6$ & 38219.5 & 125440 & 16536.1 & $3.0\times 10^{12}$\\
      \bottomrule
  \end{tabular}
\end{table} 
\noindent
Note that, for $r=2$, $k_1=2$, and $k_2=5$, the minimum range is $r(k_1+k_2)=14$. Also, we have computed the probabilities by taking $m$ from $14$ to $20$ and similarly, it can be calculated for other values of $m$.

\section{Concluding Remarks}\label{four:CR}
\begin{itemize}
\item[(i)] It is important to note that the expression $\sum_{i=1}^{m_2}(pz)^{i-1}=\sum_{i=1}^{k_2-\ell_1+1}(pz)^{i-1}$ appears in \eqref{four:ddpgf1} and \eqref{four:ddpgf3}, as expected, since the pattern can be completed if a failure occurs after $\ell_2+1$ (up to $k_2$) consecutive successes. Also, with the same justification, the expressions \eqref{four:wt2} and \eqref{four:wt3} have the term $qt/(1-pt)$. However, \eqref{four:ddpgf2} and \eqref{four:wt2} are in easy form as the pattern is completed just after $\ell_2$ consecutive successes.
\item[(ii)] The explicit form of PGF and PMF in Theorems \ref{four:ex1}, \ref{four:ex2} and \ref{four:ex3} can also be expressed in different forms as the binomial expansion can be written $(a+b)^n=\sum_{u=0}^{n}{n \choose u}a^u b^{n-u}=\sum_{u=0}^{n}{n \choose u}a^{n-u}b^u$. It is up to the end-user to choose an appropriate form and modify the results.
\item[(iii)] The results derived in Section \ref{four:DR}, are based on Markov chain approach (see Fu and Koutras \cite{FK} and Dafnis {\em et al.} \cite{DAP}). However, the results can also be derived using combinatorial method similar to Huang and Tsai \cite{HT} for i.i.d. case.
\item[(iv)] It can be easily verified that for $\ell_1=k_1$ and $\ell_2=k_2$, Theorems \ref{four:th1} - \ref{four:wth3} are same as Theorems $3.1$ - $3.8$ of Kumar and Upadhye \cite{KU3}, as expected.
\item[(v)] In Theorems $4.4$ and $4.7$, for $r \ge 1$, Dafnis {\em et al.} \cite{DAP} proved that the PGF for $r$-th waiting time of $H_{k_1,k_1,k_2,k_2}^n=X_n^{(2)}$ and $H_{1,k_1,1,k_2}^n=X_n^{(3)}$ (in their notation) are given by
%For $\ell_1=k_1$ and $\ell_2=k_2$, $H_{k_1,k_1,k_2,k_2}^n=X_n^{(2)}$ of Dafnis {et al.} \cite{DAP} (in their notation). Dafnis {et al.} \cite{DAP} in Theorem $4.4$ proved that for $r \ge 1$, the PGF for waiting time for $X_n^{(2)}$ is given by
\begin{equation}
H_r(z)=\left(\frac{(qz)^{k_1}(pz)^{k_2}(1-qz)(1-pz)}{1-z+(qz)^{k_1}(pz)^{k_2}(1-qz)(1-pz)}\right)^r(1-pz)^{-1}.\label{four:errr}
\end{equation}
%But, observe that $H_r(1)=1/(1-p) \neq 1$ unless $p=0$. Also, for $\ell_1=1=\ell_2$, $H_{1,1,k_1,k_2}^n=X_n^{(3)}$ of Dafnis {et al.} \cite{DAP} (in their notation). Also, Dafnis {et al.} \cite{DAP} in Theorem $4.7$ proved that for $r \ge 1$, the PGF for waiting time for $X_n^{(3)}$ is given by
\begin{equation}
{\rm and}~~~~~~~~~~~~~~~~~~~~~~~~~~H_r(z)=\left(\frac{(qz)(pz)(1-(qz)^{k_1})(1-(pz)^{k_2})}{1-z+(qz)(pz)(1-(qz)^{k_1})(1-(pz)^{k_2})}\right)^r(1-(pz)^{k_2})^{-1}.\label{four:err}~~~~~~~~~~~~~~~
\end{equation}
respectively. But, observe that $H_r(1)=1/(1-p) \neq 1$ in \eqref{four:errr} and $H_r(1)=1/(1-p^{k_2}) \neq 1$ in \eqref{four:err} unless $p=0$. Therefore, the expressions \eqref{four:errr} and \eqref{four:err} are incorrect and hence Theorems $4.5$, $4.6$, $4.8$ and $4.9$ of Dafnis {et al.} \cite{DAP} are also incorrect. We correct and generalize these erroneous results in Theorem \ref{four:wth3}.
% Hence, we corrected some errors of Dafnis {et al.} \cite{DAP}. Also, we can say that Theorems $3.7$ - $3.9$ may not be correct as Theorems $3.7$ - $3.9$ and Theorem $4.7$ - $4.9$ are obtained by using same matrix $A$, $B$ and $C$. However PGF and PMF of $X_n^{(3)}$, which is $\varphi_n(t)$ and $\kappa_{m,n}$ in our notation, satisfied the necessary conditions.
\item[(vi)] Note that if $p_0=1$, $p_1=0$, $p_{00}=q=p_{10}$ and $p_{01}=p=p_{11}$ then \eqref{four:Mddpgf2} $\implies$ \eqref{four:ddpgf2}, \eqref{four:Mwt2} $\implies$ \eqref{four:wt2}, \eqref{four:Mddpgf1} $\implies$ \eqref{four:ddpgf1}, \eqref{four:MMMwt1} $\implies$ \eqref{four:wt1}, \eqref{four:Mddpgf3} $\implies$ \eqref{four:ddpgf3} and \eqref{four:MMMwt3} $\implies$ \eqref{four:wt3}, as expected.
\item[(vii)] Using the double generating functions obtained in Subsection \ref{four:MDR}, the results for PGF, PMF and moments similar to  Subsections \ref{four:D1}, \ref{four:D2} and \ref{four:D3} can be derived.
\end{itemize}

%It can be concluded that if $k_1=\infty$ then {\bf (T3)} $\implies$ {\bf (T2)}. Also, this can be verified by letting $k_1=\infty$, i.e., $m_1=k_1-\ell_1+1=\infty$, in (7) which implies (5). However, we have derived explicit results for {\bf (T2)} for the sake of simplicity to the reader. But, {\bf (T3)} $\centernot\implies$ {\bf (T1)} by letting $k_2=\infty$ because 0 is needed after at least $\ell_2$ (at most $k_2$) consecutive 1’s for {\bf (T3)} but not for {\bf (T1)}.
%The results derived in Section \ref{four:MDR} are general results over Section \ref{four:DR}.

%\small
\footnotesize
\singlespacing
\section*{Acknowledgements}
The authors are grateful to the associate editor and reviewers for many valuable suggestions, critical comments which improved the presentation of the paper.

\end{document}